\documentclass[12pt,a4paper,titlepage]{article}
\usepackage{amsmath,amsfonts,amssymb,bbm,mathrsfs,theorem}
\usepackage{t1enc}
\usepackage[latin1]{inputenc}
\usepackage[english]{babel}
\usepackage{epsfig,here}
\usepackage{times}

\newcommand{\D}{\ensuremath{\mathscr{D}}}

\newcommand{\N}{\ensuremath{\mathbbm{N}}}
\newcommand{\R}{\ensuremath{\mathbbm{R}}}

\newcommand{\Rs}{\ensuremath{R^\ast}}
\newcommand{\laplace}{\ensuremath{\triangle}}

\renewcommand{\S}{\ensuremath{\mathscr{S}}}
\renewcommand{\(}{\left(}
\renewcommand{\)}{\right)}
\renewcommand{\exp}{\mbox{e}}
\newcommand{\uint}{\int\limits}

\newcommand{\ph}{\ensuremath{\varphi}}
\newcommand{\pa}{\partial}

\theoremstyle{changebreak}

\theoremheaderfont{\bf\scshape}

\newtheorem{defi}{Definition}[section]
\newtheorem{prop}[defi]{Proposition}
\newtheorem{rem}[defi]{Remark}
\newtheorem{cor}[defi]{Corollary}
\theorembodyfont{\itshape}
\newtheorem{theo}[defi]{Theorem}
\newenvironment{proof}{\noindent{\bf Proof:}\par\nopagebreak}{\hspace*{\fill} $
\Box$}

\begin{document}

\begin{sloppypar}

\renewcommand{\textfraction}{0}
\renewcommand{\topfraction}{0.9}
\renewcommand{\bottomfraction}{0.9}
\setcounter{topnumber}{10}
\setcounter{bottomnumber}{10}
\setcounter{totalnumber}{10}

{\bf \Large Inverting the spherical Radon transform for physically meaningful functions}

\begin{center}

{\large Jens Klein}\\
{\it Institut f\"ur numerische\\
und instrumentelle Mathematik\\
Westf\"alische Wilhelms-Universit\"at M\"unster\\
Einsteinstrasse 62, D-48149 M\"unster, Germany\\
e-mail: jens.klein@math.uni-muenster.de}

\end{center}

\begin{quote}
{\small {\bf Abstract}
In this paper we refer to the reconstruction formulas given in Andersson's {\it On the determination of a function from spherical averages}, which are often used in applications such as SAR\footnote{SAR: Synthetic Aperture Radar} and SONAR\footnote{SONAR: SOund NAvigation and Ranging}. We demonstrate that the first one of these formulas does not 
converge given physically reasonable assumptions. An alternative is proposed and it is shown that the second reconstruction formula is well-defined but might be difficult to compute numerically.}
\end{quote}

\section{Introduction}

The determination of a function from spherical averages is a problem often encountered in physical applications such as SAR and SONAR. The work related to this topic, which has lead to a great amount of insight and refinement today, began with the proposal of a reconstruction formula by Fawcett \cite{Fawcett}. The mathematical analysis of the problem was later improved \cite{Andersson} and two refined reconstruction formulas were derived. This sparked a host of activity that branched to several areas \cite{Ulander1}, \cite{Hellsten}, \cite{Nolan}, \cite{Borden}, \cite{Louis} so that most of today's research is based on Anderssons's ideas.\\

In \cite{Andersson} two reconstruction formulas were derived from the fourier inversion formula. But it was neglected to examine whether they are properly defined. In the following it will be shown that one of these is useless, if the function to be reconstructed has physically sensible properties, because of a divergent integral in this case. An alternative will be presented and it will be shown that the the other reconstruction formula is well defined but might me difficult to compute numerically.\\
At first, for the benefit of the reader, some results from \cite{Andersson} will be recalled. Note that besides the aforementioned problems, a few minor errors occurred, which do not essentially obscure the results in \cite{Andersson}. For a detailed analysis of these errors see \cite{Klein}.

\begin{defi}

Let $n \in \N$. Then

\begin{enumerate}

\item

$\S(\R^{n+1})$ is the Schwartz Space.

\item

$\S_e(\R^{n+1}):=\{\varphi \in \S(\R^{n+1}): \varphi(x,-y)=\varphi(x,y), \forall x \in \R^n, y \in \R \}$.

\item

$\S_r(\R^n \times \R^{n+1}):=\{\varphi \in \S(\R^{2n+1}) | \forall \mbox{ orthonormal transformations}$
$$ U:\R^{n+1} \rightarrow \R^{n+1} \forall x \in \R^n , z \in \R^{n+1}: \varphi (x,z) = \varphi (x,Uz)\} \mbox{.}$$

\item

$\S'_e(\R^{n+1})$ and $\S'_r(\R^n \times \R^{n+1})$ are the dual spaces of $\S_e(\R^{n+1})$ and $\S_r(\R^n \times \R^{n+1})$ respectively.

\item

Let $f \in \S_e(\R^{n+1})$. Then the operator $R$ is defined by
$$Rf(x,r):=\frac{1}{|S^n|}\uint_{S^n}f(x+r\xi,r\eta) \, d S_n(\xi,\eta) \mbox{.}$$

\item

$S^n$ denotes the unit sphere in $\R^{n+1}$.

\end{enumerate}

\end{defi}

For simplicity in the following sometimes only $\S$, $\S_e$, etc. is written instead of $\S(\R^{n+1})$, $\S_e(\R^{n+1})$, etc.

\begin{rem}

It is easily seen that $\hat{f} \in \S_e$ and $\hat{g} \in \S_r$, if $f \in \S_e$ and $g \in \S_r$, respectively.

\end{rem}

The essential result in \cite{Andersson} is the (Fourier-) inversion formula:

\begin{theo}

\label{inversionformula}

If $\S_e(\R^{n+1})$ is given the topology of $\S'(\R^{n+1})$, the mapping
$$R:\S_e(\R^{n+1}) \rightarrow \S'_r(\R^n \times \R^{n+1})$$
is continuous and can, by continuity, be extended to a mapping
$$R:\S'_e(\R^{n+1}) \rightarrow \S'_r(\R^n \times \R^{n+1}) \mbox{.}$$
The range of this extended mapping $R$ is the closed subspace
$$\S'_{r,cone}(\R^n \times \R^{n+1})=\{g \in \S'_r(\R^n \times \R^{n+1}): \mbox{supp }\hat{g} \subseteq \{(\xi,\eta):\|\eta\| \geq \|\xi\|\}\} \mbox{.}$$
$R$ is one-to-one and the inverse mapping
$$R^{-1}:\S'_{r,cone}(\R^n \times \R^{n+1}) \rightarrow \S'_e(\R^{n+1})$$
is continuous. Moreover, if $g=Rf$ and if $\hat{f}(\xi,\eta)$ or $\hat{g}(\xi,\eta)$ are integrable for $\xi, \eta$ in some open set, then
$$\hat{g}(\xi,\eta)=\left\{ \begin{array}{lr}
\displaystyle{(2\pi)^n \frac{2}{|S^n|} \frac{\hat{f}(\xi,\sqrt{\|\eta\|^2-\|\xi\|^2})}{\|\eta\|^{n-1} \sqrt{\|\eta\|^2-\|\xi\|^2}}} &
\mbox{for }\|\eta\|>\|\xi\| ,\\
0 & \mbox{for }0 < \|\eta\| \leq \|\xi\| \\
\end{array} \right.$$
or
$$\hat{f}(\xi,\eta)=\frac{1}{(2\pi)^n} \frac{|S^n|}{2}|\eta|(\|\xi\|^2+\eta^2)^{\frac{n-1}{2}}\hat{g}(\xi,\sqrt{\|\xi\|^2+\eta^2}) \mbox{,}$$
respectively.

\end{theo}

\begin{proof}

\cite[Theorem 2.1]{Andersson} 

\end{proof}

In this theorem the operator $R$ describes in the case of SAR and SONAR the measurement of the reflectivity function $f$ that represents the ground reflectivity. The measurement is modelled as a $\delta$-impulse wave that propagates as concentric spheres. The ground is approximated as a plane. The single scatter approximation for the wave hitting the ground results in integrals over circles. The Fourier transform $\hat{g}$ of the data $g$ can under certain conditions be used to extract the Fourier-transform $\hat{f}$ of the reflectivity function $f$.

\begin{rem}

Let $f \in \S_e(\R^{n+1})$. Then for $g=Rf \in \S'_r(\R^n \times \R^{n+1})$, $x \in \R^n$ and $r \geq 0$ in the following $g(x,r)$ is sometimes written with abuse of notation. This is justified because $g$ depends only radially on the last $n+1$ variables.

\end{rem}

\begin{defi}

Let $g \in \S_r(\R^n \times \R^{n+1})$. Then 
$$\Rs g(x,y)=\uint_{\R^n} g(z,\sqrt{\|z-x\|^2+y^2}) \, dz \mbox{.}$$

\end{defi}

In \cite{Andersson} two reformulations of this formula were proposed, of which the first was already essentially given in \cite{Fawcett}.

\begin{cor}

With $c_n=\frac{1}{(2\pi)^n}\frac{|S^n|}{2}$ two reformulations of the inversion formula are possible:

\begin{enumerate}

\item

For $g \in \S_r$
$$f=c_n H_y \frac{\pa}{\pa y} \laplace^\frac{n-1}{2} \Rs g$$
with the Hilbert transform in $y$, $H_y$, and the Laplace-Operator $\laplace=\laplace_x+\frac{\pa^2}{\pa y^2}$. This formula is essentially also given by Fawcett \cite{Fawcett}.

\item

For $g \in \S_{r,cone}:=\{g \in \S_r: \mbox{supp }\hat{g} \subseteq \{(\xi,\eta):\|\eta\| \geq \|\xi\|\}\}$
$$f=c_n \Rs K g$$
with the pseudodifferential operator $K$ defined by $\widehat{Kg}(\xi,\eta)=\sqrt{\|\eta\|^2-\|\xi\|^2}\,\|\eta\|^{n-1}\hat{g}(\xi,\eta)$.

\end{enumerate}

\end{cor}

\begin{proof}

\cite[Section 3]{Andersson}

\end{proof}

The corollary states that under certain restrictions it is possible to reconstruct the reflectivity function $f$ directly from the data $g$ without taking the detour through the Fourier space.

\begin{rem}

Note that the essential restriction of this formulation of the corollary in comparison to \cite{Andersson} is that the data $g$ has to be in $\S_r$ or $S_{r, cone}$ respectively. This is necessary, because otherwise the application of $\Rs$ to $g$ in the first case is undefined or the application of $\Rs$ to $Kg$ in the second case.\\
Unfortunately $g=Rf$ is usually not in $\S_r$. Therefore the two reformulations of the inversion formula are only valid in the distributional sense with an appropriately defined $\Rs$.\\
For example a physically reasonable $f \in \S_e$, $f:\R^2 \rightarrow [0, \infty)$ with $f(x_0,y_0)>0$ for some $x_0,y_0 \in \R$
yields $f(x,y)>c>0$ for all $(x,y) \in K_\epsilon (x_0,y_0)$ with appropriate $c,\epsilon >0$.\\
Moreover
$$2\pi g(z,\sqrt{(z-x)^2+y^2})=\frac{1}{\sqrt{(z-x)^2+y^2}} \uint_{\|r\|=\sqrt{(z-x)^2+y^2}} f(z+r_1,r_2) \, d\sigma(r)$$
with $r=(r_1,r_2)$.
For $\sqrt{(z-x_0)^2+y_0^2}>\epsilon$ we obtain with a simple geometrical consideration:
$$2\pi g(z,\sqrt{(z-x_0)^2+y_0^2}) \geq \frac{\epsilon c}{2 \sqrt{(z-x_0)^2+y_0^2}} \mbox{.}$$
$$\Rightarrow 2\pi \Rs g(x_0,y_0)=\uint_\R g \( z,\sqrt{(z-x_0)^2+y_0^2} \) \, dz$$
$$\geq \frac{\epsilon c}{2} \uint_{\sqrt{(z-x_0)^2+y_0^2}>\epsilon}
\frac{1}{\sqrt{(z-x_0)^2+y_0^2}} \, dz =\infty \mbox{.}$$
Therefore the reconstruction of a non-negative function $f$ with $f(x_0,y_0)>0$ for some $(x_0,y_0) \in \R^2$ is impossible using the first reconstruction formula.

\end{rem}

This result is in accordance with a result from Nessibi, Rachdi and Trimeche \cite{Nessibi}. They gave reconstruction formulas for functions $f$ with
$$\uint_0^\infty P(r) f(r,x) \, dr=0$$
for all $x \in \R^n$ and for all one-variable polynomials $P$.

\section{Properties of the function $g=R f$}

Before showing an important property of the data function $g=Rf$, which will be necessary for the derivation of the new reconstruction formula, some definitions and a result from \cite{Andersson} will be needed.

\begin{defi}

Let $n \in \N$. Then

\begin{enumerate}

\item

$\D(\R^{n+1}):=C_0^\infty(\R^{n+1})$.

\item

$\D_e(\R^{n+1}):=\{\varphi \in \D(\R^{n+1}): \varphi(x,-y)=\varphi(x,y), \forall x \in \R^n, y \in \R \}$.

\end{enumerate}

\end{defi}

\begin{cor}

\label{dipl4.3.2}

If $f \in \S_e(\R^{n+1})$, then $g=Rf \in \S'_r(\R^n \times \R^{n+1})$.

\end{cor}

\begin{proof}

\cite[Section 2]{Andersson}

\end{proof}

\begin{prop}

\label{gincinf}

If $f \in \S_e(\R^{n+1})$, then $g=Rf \in C^\infty$.

\end{prop}

\begin{proof}

The change of differentiation and integration is justified, because $f \in \S_e(\R^{n+1})$.

\end{proof}

\begin{theo}

\label{ghatinl1}

If $f \in \S_e(\R^{n+1})$, then $\hat{g}=\widehat{Rf} \in L^1(\R^n \times \R^{n+1})$.

\end{theo}

\begin{proof}

With $f \in \S$, also $\hat{f} \in \S \subseteq L^1$. Theorem \ref{inversionformula} implies that
$$\hat{g}(\xi,\eta)=\left\{ \begin{array}{lr}
\displaystyle{(2\pi)^n \frac{2}{|S^n|} \frac{\hat{f}(\xi,\sqrt{\|\eta\|^2-\|\xi\|^2})}{\|\eta\|^{n-1} \sqrt{\|\eta\|^2-\|\xi\|^2}}} &
\mbox{for }\|\eta\|>\|\xi\| ,\\
0 & \mbox{otherwise.} \\
\end{array} \right.$$
$$\Rightarrow \uint_{\R^n \times \R^{n+1}} |\hat{g}(\xi,\eta)|\,d\xi d\eta=(2\pi)^n\frac{2}{|S^n|}\uint_{\|\eta\| \geq \|\xi\|}\left|\frac{\hat{f}(\xi,\sqrt{\|\eta\|^2-\|\xi\|^2})}{\|\eta\|^{n-1} \sqrt{\|\eta\|^2-\|\xi\|^2}}\right|\,d\xi d\eta \mbox{.}$$
The substitution $\rho'=\|\eta\|$ results in
$$=2(2\pi)^n\uint_{\R^{n}}\uint_{\rho' \geq \|\xi\|}\left|\frac{\rho' \hat{f}(\xi,\sqrt{\rho'^2-\|\xi\|^2})}{\sqrt{\rho'^2-\|\xi\|^2}}\right|\,d\xi d\rho' \mbox{.}$$
The substitution $\rho=\rho' + \|\xi\|$ leads to
$$=2(2\pi)^n\uint_{\R^n}\uint_{\rho \geq 0}\left|\frac{(\rho + \|\xi\|) \hat{f}(\xi,\sqrt{\rho^2+2\rho \|\xi\|})}{\sqrt{\rho^2+2\rho \|\xi\|}}\right|\,d\xi d\rho \mbox{.}$$
$f$ is in $\S$, therefore
$$\leq 2(2\pi)^n\uint_{\R^n}\uint_{\rho \geq 0}\frac{C (\rho + \|\xi\|)}{\sqrt{\rho^2+2\rho \|\xi\|} (1+\sqrt{\|\xi\|^2+\rho^2+2\rho \|\xi\|})^{n+3}}\,d\xi d\rho$$
$$=2(2\pi)^n\uint_{\R^n}\uint_{\rho \geq 0}\frac{C (\rho+\|\xi\|)}{\sqrt{\rho^2 + 2\rho \|\xi\|}(1+\|\xi\|+\rho)^{n+3}} \,d\xi d\rho$$
$$\leq 2(2\pi)^n\uint_{\R^n}\uint_{\rho \geq 0} \frac{C}{\sqrt{\rho^2 + 2\rho \|\xi\|}(1+\|\xi\|+\rho)^{n+2}} \,d\xi d\rho$$
and the integral
$$\uint_{\R^n}\uint_{\rho \geq 1} \frac{C}{\sqrt{\rho^2 + 2\rho \|\xi\|}(1+\|\xi\|+\rho)^{n+2}} \,d\xi d\rho$$
converges.
So only the integral
$$\uint_{\R^n}\uint_{0 \leq \rho \leq 1} \frac{C}{\sqrt{\rho^2 + 2\rho \|\xi\|}(1+\|\xi\|+\rho)^{n+2}} \,d\xi d\rho$$
remains to be examined.
$$\uint_{\R^n}\uint_{0 \leq \rho \leq 1} \frac{C}{\sqrt{\rho^2 + 2\rho \|\xi\|}(1+\|\xi\|+\rho)^{n+2}} \,d\xi d\rho$$
$$\leq \uint_{\R^n}\uint_{0 \leq \rho \leq 1} \frac{C}{\sqrt{\rho \|\xi\|}(1+\|\xi\|)^{n+2}} \,d\xi d\rho$$
$$=\uint_{\R^n}\frac{2C}{\sqrt{\|\xi\|}(1+\|\xi\|)^{n+2}} \,d\xi \mbox{.}$$
The integral
$$\uint_{\|\xi\| \geq 1}\frac{2C}{\sqrt{\|\xi\|}(1+\|\xi\|)^{n+2}} \,d\xi$$
converges. For the proof it therefore suffices to show the existence of
$\uint_{\|\xi\| \leq 1}\frac{2C}{\sqrt{\|\xi\|}(1+\|\xi\|)^{n+2}} \,d\xi \mbox{:}$
$$\uint_{\|\xi\| \leq 1}\frac{2C}{\sqrt{\|\xi\|}(1+\|\xi\|)^{n+2}} \,d\xi$$
$$\leq \uint_{\|\xi\| \leq 1} \frac{2C}{\sqrt{\|\xi\|}} \, d\xi \mbox{.}$$
The substitution $r=\|\xi\|$ yields
$$=2C \left| S^{n-1} \right| \uint_0^1 r^{n-\frac{3}{2}} \, dr < \infty \mbox{,}$$
because $n \geq 1$.

\end{proof}

\begin{rem}

An analogous proof shows $\eta \hat{g} \in L^1$.

\end{rem} 

\section{Modifications to Andersson's first inversion formula}

\subsection{Definition and properties of a modified $\Rs$}

Now a modified version of the operator $\Rs$ is introduced and its properties are discussed.

\begin{defi}

\label{rsternpa}

For $f \in \D_e$ and $g=Rf \in C^\infty(\R^n \times \R^{n+1})$ we define
$$(\Rs_{\pa} g)(x,y):=\uint_{\R^n} \frac{\pa}{\pa y} g(z,\sqrt{\|x-z\|^2 + y^2})\, dz \mbox{.}$$

\end{defi} 

This slight modification by an additional derivation turns out to ensure the convergence of the integral applied by the operator $\Rs_{\pa}$ under the minor and physically feasible constraint that $f$ is in $\D_e$. Thereby the formulation of a mathematically exact reconstruction formula is possible.

\begin{cor}

\label{rsternpawohldef}

$\Rs_{\pa} g$ is well-defined.

\end{cor}

\begin{proof}

Let $x \in \R^n$, $y \in \R$. Then
$$\left|(\Rs_{\pa} g)(x,y)\right|=\left|\uint_{\R^n} \frac{\pa}{\pa y} g(z,\sqrt{\|x-z\|^2 + y^2})\, dz\right|$$
$$=\left|\uint_{\R^n} \frac{\pa}{\pa y} \frac{1}{|S^n|} \uint_{S^n} f(z+\xi\sqrt{\|x-z\|^2 + y^2},\eta\sqrt{\|x-z\|^2 + y^2})\, dS_n(\xi,\eta)dz\right|\mbox{.}$$
$f \in \D_e$ and therefore it is possible to exchange differentiation and integration.
$$=\left|\uint_{\R^n} \frac{1}{|S^n|} \uint_{S^n} \frac{\pa}{\pa y} f(z+\xi\sqrt{\|x-z\|^2 + y^2},\eta\sqrt{\|x-z\|^2 + y^2})\, dS_n(\xi,\eta)dz\right|$$
$$=\left|\uint_{\R^n} \frac{1}{|S^n|} \uint_{S^n} \frac{y}{\sqrt{\|x-z\|^2 + y^2}}{\xi \choose \eta} \cdot (\nabla f)(z+\xi\sqrt{\|x-z\|^2 + y^2},\eta\sqrt{\|x-z\|^2 + y^2})\, dS_n(\xi,\eta)dz\right|$$
$$\leq \uint_{\R^n} \frac{1}{|S^n|} \uint_{S^n} \frac{|y|}{\sqrt{\|x-z\|^2 + y^2}} \|(\nabla f)(z+\xi\sqrt{\|x-z\|^2 + y^2},\eta\sqrt{\|x-z\|^2 + y^2})\|\, dS_n(\xi,\eta)dz \mbox{.}$$
With $r=(\xi,\eta)\sqrt{\|x-z\|^2 + y^2}$ this is
$$=\uint_{\R^n} \frac{1}{|S^n|} \uint_{\|r\|=\sqrt{\|x-z\|^2 + y^2}} \frac{|y|}{(\sqrt{\|x-z\|^2 + y^2}\,)^{n+1}} \|(\nabla f)((z,0)+r)\|\, d\sigma(r)dz \mbox{.}$$
$f \in \D \Rightarrow |(\Rs_{\pa} g)(x,y)|$
$$\leq \uint_{\R^n} \frac{|y|}{(\sqrt{\|x-z\|^2 + y^2}\,)^{n+1}} \mbox{ max }(\|\nabla f\|) \mbox{ diam}(\mbox{supp } f) \, dz < \infty$$

\end{proof}

\begin{cor}

\label{rsternpainsstrich}

For $f \in \D_e$ and $g=Rf$, $\Rs_{\pa}g \in \S'$.

\end{cor}

\begin{proof}

This is guaranteed by the estimate in Corollary \ref{rsternpawohldef}.

\end{proof}

\begin{defi}

Let $f \in \S'(\R^n)$ and $\ph \in \S(\R^n)$. Then
$$<f,\ph>_{\S(\R^n)}:=f(\ph)$$
is the functional $f$ applied to the test function $\ph$.
Here the subscript $\S(\R^n)$ is a reminder of the space of the test function $\ph$.

\end{defi}

Now analogously to \cite{Andersson} an expression for $\widehat{\Rs_{\pa} g}$ is derived. 

\begin{theo}

\label{rsternpagdachgleichetagdach}

Let $f \in \D_e$ and $g=Rf$. Then $\widehat{\Rs_{\pa}g}(\xi,\eta)=i \eta \hat{g}(\xi,\sqrt{\|\xi\|^2+\eta^2})$.

\end{theo}

\begin{proof}

Let $\ph\in \S$, $f \in \D_e$ and $g=Rf$. Then
$$<\widehat{\Rs_{\pa}g},\ph>_{\S(\R^n \times \R)}=<\Rs_{\pa}g,\hat{\ph}>_{\S(\R^n \times \R)}=\uint_{\R^n \times \R} (\Rs_{\pa}g)(x',y)\hat{\ph}(x',y)\,dx'dy$$
$$=\uint_{\R^n \times \R} \uint_{\R^n} \frac{\pa}{\pa y} g(z,\sqrt{\|x'-z\|^2 + y^2})\, dz \hat{\ph}(x',y)\,dx'dy \mbox{.}$$
As $\Rs_{\pa}g \in \S'$ and $\hat{\ph} \in \S$, Fubini's theorem implies
$$=\uint_{\R^n} \uint_{\R^n \times \R} \frac{\pa}{\pa y} g(z,\sqrt{\|x'-z\|^2 + y^2}) \hat{\ph}(x',y)\,dx'dydz \mbox{.}$$
Substituting $x=x'+z$ we obtain
$$=\uint_{\R^n} \uint_{\R^n \times \R} \frac{\pa}{\pa y} g(z,\sqrt{\|x\|^2 + y^2}) \hat{\ph}(x+z,y)\,dxdydz \mbox{.}$$
Let $\psi(x,y,z):=\hat{\ph}(x+z,y)$. Then $\psi(x,y,z) \in \S(\R^n \times \R \times \R^n)$ and this leads to
$$=\uint_{\R^n} \uint_{\R^n \times \R} \frac{\pa}{\pa y} g(z,\sqrt{\|x\|^2 + y^2}) \psi(x,y,z)\,dxdydz$$
$$=<\frac{\pa}{\pa y}g,\psi>_{\S(\R^n \times \R \times \R^n)}=<\widehat{\frac{\pa}{\pa y}g},\tilde{\psi}>_{\S(\R^n \times \R \times \R^n)}=<i \eta\hat{g},\tilde{\psi}>_{\S(\R^n \times \R \times \R^n)}$$
$$=\uint_{\R^n \times \R} \uint_{\R^n} i \eta \hat{g}(\zeta,\sqrt{\|\xi\|^2+\eta^2})\tilde{\psi}(\xi,\eta,\zeta) \, d\zeta d\xi d\eta$$
$$=\uint_{\R^n \times \R} \uint_{\R^n} i \eta \hat{g}(\zeta,\sqrt{\|\xi\|^2+\eta^2}) \uint_{\R^n} \uint_{\R^n \times \R} (2\pi)^{-2n-1} \exp^{i(<x,\xi>+y\eta+<z,\zeta>)} \psi(x,y,z) \, dxdydz d\zeta d\xi d\eta$$
$$=\uint_{\R^n \times \R} \uint_{\R^n} i \eta \hat{g}(\zeta,\sqrt{\|\xi\|^2+\eta^2}) \uint_{\R^n} \uint_{\R^n \times \R} (2\pi)^{-2n-1} \exp^{i(<x,\xi>+y\eta+<z,\zeta>)} \hat{\ph}(x+z,y) \, dxdydz d\zeta d\xi d\eta \mbox{.}$$
As $\eta\hat{g} \in L^1$ and $\hat{\ph} \in \S$, Fubini's theorem gives
$$=\uint_{\R^n} \uint_{\R^n \times \R} \uint_{\R^n \times \R} \uint_{\R^n} i \eta \hat{g}(\zeta,\sqrt{\|\xi\|^2+\eta^2})(2\pi)^{-2n-1} \exp^{i(<x,\xi>+y\eta+<z,\zeta>)} \hat{\ph}(x+z,y) \, d\zeta d\xi d\eta dxdydz \mbox{.}$$
The substitution $x'=x-z$ yields
$$=\uint_{\R^n} \uint_{\R^n \times \R} \uint_{\R^n \times \R} \uint_{\R^n} i \eta \hat{g}(\zeta,\sqrt{\|\xi\|^2+\eta^2}) (2\pi)^{-2n-1} \exp^{i(<x',\xi>+y\eta+<z,\zeta-\xi>)} \hat{\ph}(x',y) \, d\zeta d\xi d\eta dx'dydz$$
$$=\uint_{\R^n} \uint_{\R^n \times \R} \uint_{\R^n} i \eta \hat{g}(\zeta,\sqrt{\|\xi\|^2+\eta^2}) (2\pi)^{-n} \exp^{i<z,\zeta-\xi>} \ph(\xi,\eta) \, d\zeta d\xi d\eta dz \mbox{.}$$
Continuing in the distributional sense we obtain
$$=\uint_{\R^n \times \R} \uint_{\R^n} i \eta \hat{g}(\zeta,\sqrt{\|\xi\|^2+\eta^2}) \delta(\zeta-\xi) \ph(\xi,\eta) \, d\zeta d\xi d\eta$$
$$=\uint_{\R^n \times \R} i \eta \hat{g}(\xi,\sqrt{\|\xi\|^2+\eta^2}) \ph(\xi,\eta) \, d\xi d\eta$$
$$=<i \eta \hat{g_1},\ph>_{\S(\R^n \times \R)}\mbox{.}$$
Here $\hat{g_1}(\xi,\eta):=\hat{g}(\xi,\sqrt{\|\xi\|^2+\eta^2})$ and with this definition $\hat{g_1} \in \S'(\R^n \times \R)$.

\end{proof}

\subsection{A modified inversion formula}

With the results of the preceding sections a well-defined reconstruction formula is attainable also for physically meaningful reflectivity functions.

\begin{theo}

\label{modinvform}

Let $f \in \D_e$ and $g=Rf$. Then $f=c_n H_y \laplace^{\frac{n-1}{2}} \Rs_{\pa} g$ with the constant $c_n:=\frac{1}{(2\pi)^n}\frac{|S^n|}{2}$ and the Hilbert-transform $H_y$.

\end{theo}

\begin{proof}

Let $f \in \D_e$. It follows from Theorem \ref{ghatinl1} and Theorem \ref{inversionformula} that
$$f=\tilde{\hat{f}}=c_n(|\eta|(\|\xi\|^2+\eta^2)^{\frac{n-1}{2}}\hat{g})\tilde{}$$
$$=c_n (-i \, sgn(\eta)(\|\xi\|^2+\eta^2)^{\frac{n-1}{2}} i \eta \hat{g})\tilde{} \mbox{.}$$
Theorem \ref{rsternpagdachgleichetagdach} yields
$$=c_n (-i \, sgn(\eta)(\|\xi\|^2+\eta^2)^{\frac{n-1}{2}} \widehat{\Rs_{\pa} g}) \tilde{}$$
$$=c_n (-i \, sgn(\eta) (\laplace^{\frac{n-1}{2}} \Rs_{\pa} g)\hat{}) \tilde{}$$
$$=c_n ((H_y \laplace^{\frac{n-1}{2}} \Rs_{\pa} g)\hat{}) \tilde{}$$
$$=c_n H_y \laplace^{\frac{n-1}{2}} \Rs_{\pa} g \mbox{.}$$

\end{proof} 

\section{Well-definedness of Andersson's second inversion formula}

In the following it is shown that the second reformulation of the inversion formula in \cite{Andersson} is well-defined. To this end at first some properties of $\widehat{Kg}$ are derived.

\begin{prop}

\label{kgdachinlp}

Let $f \in \S_e$ and $g=Rf$. Then $\widehat{Kg} \in L^p(\R^n \times \R^{n+1})$ for all $p>0$. In general $g$ is only in $C^\infty$ and not for example in $L^2$. Therefore $\hat{g}$ has to be computed in the distributional sense.

\end{prop}

\begin{proof}

Theorem \ref{inversionformula} implies that
$$\hat{g}(\xi,\eta)=\left\{ \begin{array}{lr}
\displaystyle{(2\pi)^n \frac{2}{|S^n|} \frac{\hat{f}(\xi,\sqrt{\|\eta\|^2-\|\xi\|^2})}{\|\eta\|^{n-1} \sqrt{\|\eta\|^2-\|\xi\|^2}}} &
\mbox{for }\|\eta\|>\|\xi\| ,\\
0 & \mbox{otherwise} \\
\end{array} \right.$$
$$\Rightarrow \left| \widehat{Kg}(\xi,\eta)\right|=\left\{ \begin{array}{lr}
\displaystyle{(2\pi)^n \frac{2}{|S^n|} \hat{f}(\xi,\sqrt{\|\eta\|^2-\|\xi\|^2})} & \mbox{for }\|\eta\|>\|\xi\| ,\\
0 & \mbox{otherwise} \\
\end{array} \right.$$
$$\leq \left\{ \begin{array}{lr}
\displaystyle{(2\pi)^n \frac{2}{|S^n|} (1+\sqrt{\|\xi\|^2+|\|\eta\|^2-\|\xi\|^2|})^{-\frac{2n+2}{p}}} & \mbox{for }\|\eta\|>\|\xi\| ,\\
0 & \mbox{otherwise} \\
\end{array} \right.$$
$$\leq \left\{ \begin{array}{lr}
\displaystyle{(2\pi)^n \frac{2}{|S^n|} (1+\|\eta\|)^{-\frac{2n+2}{p}}} & \mbox{for }\|\eta\|>\|\xi\| ,\\
0 & \mbox{otherwise} \\
\end{array} \right.$$
$$\leq (2\pi)^n \frac{2}{|S^n|} (1+\mbox{max}(\|\xi\|,\|\eta\|))^{-\frac{2n+2}{p}} \mbox{.}$$

\end{proof}

\begin{cor}

Let $f \in \D_e$ with $\hat{f}(0,0) \not= 0$ and $g=Rf$. Then $g,\hat{g} \not\in L^2(\R^n \times \R^{n+1})$.

\end{cor}

\begin{proof}

Theorem \ref{inversionformula} yields
$$\hat{g}(\xi,\eta)=\left\{ \begin{array}{lr}
\displaystyle{(2\pi)^n \frac{2}{|S^n|} \frac{\hat{f}(\xi,\sqrt{\|\eta\|^2-\|\xi\|^2})}{\|\eta\|^{n-1} \sqrt{\|\eta\|^2-\|\xi\|^2}}} &
\mbox{for }\|\eta\|>\|\xi\| ,\\
0 & \mbox{otherwise} \\
\end{array} \right.$$
$$\Rightarrow \uint_{\R^n \times \R^{n+1}} |\hat{g}(\xi,\eta)|^2\,d\xi d\eta=(2\pi)^{2n}\frac{4}{|S^n|^2}\uint_{\|\eta\| \geq \|\xi\|}\frac{|\hat{f}(\xi,\sqrt{\|\eta\|^2-\|\xi\|^2})|^2}{\|\eta\|^{2n-2} (\|\eta\|^2-\|\xi\|^2)}\,d\xi d\eta \mbox{.}$$
The substitution $\|\eta\|=\rho'$ results in
$$=(2\pi)^{2n}\frac{4}{|S^n|}\uint_{\R^{n}}\uint_{\rho' \geq \|\xi\|}\frac{\rho'^{2-n}|\hat{f}(\xi,\sqrt{\rho'^2-\|\xi\|^2})|^2}{\rho'^2-\|\xi\|^2}\,d\xi d\rho' \mbox{.}$$
The substitution $\rho'=\rho + \|\xi\|$ leads to
$$=(2\pi)^{2n}\frac{4}{|S^n|}\uint_{\R^n}\uint_{\rho \geq 0}\frac{(\rho + \|\xi\|)^{2-n} |\hat{f}(\xi,\sqrt{\rho^2+2\rho \|\xi\|})|^2}{\rho^2+2\rho \|\xi\|}\,d\xi d\rho \mbox{.}$$
We consider
$$\uint_{\|\xi\| \leq \frac{1}{2}}\uint_{0 \leq \rho \leq \frac{1}{2}}\frac{(\rho + \|\xi\|)^{2-n} |\hat{f}(\xi,\sqrt{\rho^2+2\rho \|\xi\|})|^2}{\rho^2+2\rho \|\xi\|}\,d\xi d\rho$$
and assume without loss of generality $\hat{f}(\xi,\sqrt{\rho^2+2\rho\|\xi\|}) \geq 1$ for $\|\xi\| \leq \frac{1}{2}, 0 \leq \rho \leq \frac{1}{2}$:
$$\geq \uint_{\|\xi\| \leq \frac{1}{2}}\uint_{0 \leq \rho \leq \frac{1}{2}}\frac{(\rho + \|\xi\|)^{2-n} }{\rho^2+2\rho \|\xi\|}\,d\xi d\rho \mbox{.}$$
With the substitution $r=\|\xi\|$ we obtain
$$=|S^{n-1}|\uint_{0 \leq r \leq \frac{1}{2}}\uint_{0 \leq \rho \leq \frac{1}{2}}\frac{(\rho + r)^{2-n} r^{n-1}}{\rho^2+2\rho r}\,dr d\rho \mbox{.}$$
For $n=1$ this is
$$\geq \uint_{0 \leq r \leq \frac{1}{2}}\uint_{0 \leq \rho \leq \frac{1}{2}}\frac{\rho + r}{\rho^2+2\rho r}\,dr d\rho
\geq \uint_{0 \leq r \leq \frac{1}{2}}\uint_{0 \leq \rho \leq \frac{1}{2}}\frac{1}{2\rho}\,dr d\rho=\infty \mbox{.}$$
For $n \geq 2$ this results in
$$\geq \uint_{0 \leq r \leq \frac{1}{2}}\uint_{0 \leq \rho \leq \frac{1}{2}}\frac{r^{n-1}}{\rho^2+2\rho r}\,dr d\rho=\infty \mbox{.}$$

\end{proof}

\begin{rem}

Unfortunately $\hat{f}(0,0) \not= 0$ for a reflectivity function $f$ with physically feasible properties. Moreover $\hat{g}$ is not continuous and therefore $g \notin L^1$, so it might be difficult to compute $\hat{g}$ numerically with sufficient accuracy.

\end{rem}

\begin{theo}

\label{rsternkgdachwelldef}

Let $f \in \S_e$ and $g=Rf$. Then $\widehat{\Rs Kg}$ is well-defined with $\hat{g}$ computed in the distributional sense. Further
$$\widehat{\Rs Kg}(\xi,\eta)=\widehat{Kg}(\xi,\sqrt{\|\xi\|^2+\eta^2}) \mbox{.}$$

\end{theo}

\begin{proof}

It follows from Proposition \ref{kgdachinlp} that $\widehat{Kg} \in L^1 \cap L^2$. Therefore
$$\widehat{Kg}(\xi,\sqrt{\|\xi\|^2+\eta^2})$$
$$=\uint_{\R^n} \( \exp^{-i<\xi,z>} \uint_{\R^n \times \R} \exp^{-i(<\xi,x>+\eta y)} Kg(z,\sqrt{\|x\|^2+y^2})\,dxdy \) dz$$
$$=\uint_{\R^n} \( \exp^{-i<\xi,z>} \uint_{\R^n \times \R} \exp^{-i(<\xi,x'-z>+\eta y)} Kg(z,\sqrt{\|x'-z\|^2+y^2})\,dx'dy \) dz$$
$$=\uint_{\R^n \times \R} \uint_{\R^n} \exp^{-i<\xi,z>} \exp^{-i(<\xi,x'-z>+\eta y)} Kg(z,\sqrt{\|x'-z\|^2+y^2})\,dzdx'dy$$
$$=\uint_{\R^n \times \R} \exp^{-i(<\xi,x'>+\eta y)} \underbrace{\uint_{\R^n} Kg(z,\sqrt{\|x'-z\|^2+y^2})\,dz}_{=\Rs Kg(x',y)}\, dx'dy$$
$$=\widehat{\Rs Kg}(\xi,\eta)$$
and the change in the order of integration and the existence of $\widehat{\Rs Kg}$ is justified by Fubini's Theorem.

\end{proof}

With this result it can be shown that the second reformulation of the inversion formula in \cite{Andersson} is well-defined.

\begin{cor}

\label{secondinversionformula}

Let $f \in \S_e$ and $g=Rf$. Then $f(x,y)=c_n \Rs Kg(x,y)$ with a constant $c_n$.

\end{cor}

\begin{proof}

Theorem \ref{inversionformula} implies that 
$$\hat{f}(\xi,\eta)=c_n |\eta|(\|\xi\|^2+\eta^2)^\frac{n-1}{2}\hat{g}(\xi,\sqrt{\|\xi\|^2+\eta^2})$$
with a constant $c_n$.
$$\Rightarrow \hat{f}(\xi,\eta)=c_n \widehat{Kg}(\xi,\sqrt{\|\xi\|^2+\eta^2})$$
$$=c_n \widehat{\Rs Kg}(\xi,\eta) \mbox{.}$$
This completes the proof, because $f$ is in $S_e$.

\end{proof} 

\onecolumn

\end{sloppypar}

\end{document}